\documentclass[runningheads]{llncs}
\usepackage{amsmath,amsfonts,amssymb}
\usepackage{graphicx}
\begin{document}
\title{The Riemannian barycentre as a \\ proxy for global optimisation}

\author{Salem Said\inst{1} \and
Jonathan H. Manton\inst{2}}

\institute{Laboratoire IMS (CNRS 5218), Universit\'e de Bordeaux 
\email{salem.said@u-bordeaux.fr}\\
\and
Department of Electrical and Electronic Engineering,\\ The University of Melbourne\\
\email{j.manton@ieee.org}}

\authorrunning{S. Said \& J. H. Manton}
\titlerunning{Riemannian barycentres for global optimisation}
\maketitle              
\begin{abstract} Let $M$ be a simply-connected compact Riemannian symmetric space, and $U$ a twice-differentiable function on $M$, with unique global minimum at $x^* \in M$. The idea of the present work is to replace the problem of searching for the global minimum of $U$, by the problem of finding the Riemannian barycentre of the Gibbs distribution $P_{\scriptscriptstyle{T}} \propto \exp(-U/T)$. In other words, instead of minimising the function $U$ itself, to minimise $\mathcal{E}_{\scriptscriptstyle{T}}(x) = \frac{1}{2}\int d^{\scriptscriptstyle 2}(x,z)P_{\scriptscriptstyle{T}}(dz)$, where $d(\cdot,\cdot)$ denotes Riemannian distance. The following original result is proved\,: if $U$ is invariant by geodesic symmetry about $x^*$, then for each $\delta < \frac{1}{2} r_{\scriptscriptstyle cx}$ ($r_{\scriptscriptstyle cx}$ the convexity radius of $M$),\\ there exists $T_{\scriptscriptstyle \delta}$ such that $T \leq T_{\scriptscriptstyle \delta}$ implies $\mathcal{E}_{\scriptscriptstyle{T}}$ is strongly convex on the geodesic ball $B(x^*,\delta)\,$, and $x^*$ is the unique global minimum of $\mathcal{E}_{\scriptscriptstyle{T\,}}$. Moreover, this $T_{\scriptscriptstyle \delta}$ can be computed explicitly. This result gives rise to a general algorithm for black-box optimisation, which is briefly described, and will be further explored in future work. 
\keywords{\hspace*{-0.1cm}Riemannian barycentre\,$\cdot$\,black-box optimisation\,$\cdot$\,symmetric space.}
\end{abstract}
\vspace{-0.5cm}
It is common knowledge that the Riemannian barycentre $\bar{x}$, of a probability distribution $P$ defined on a Riemannian manifold $M$, may fail to be unique. However, if $P$ is supported inside a geodesic ball $B(x^*,\delta)$ with radius $\delta < \frac{1}{2} r_{\scriptscriptstyle cx}$ ($r_{\scriptscriptstyle cx}$ the convexity radius of $M$), then $\bar{x}$ is unique and also belongs to  $B(x^*,\delta)$. In~fact,~Afsari has shown this to be true, even when $\delta < r_{\scriptscriptstyle cx}$ (see~\cite{karcher}\cite{afsari}). 

\indent Does this statement continue to hold, if $P$ is not supported inside $B(x^*,\delta)$, but merely concentrated on this ball? The answer to this question is positive, assuming that $M$ is a simply-connected compact Riemannian symmetric space, and $P = P_{\scriptscriptstyle{T}} \propto \exp(-U/T)$, where the function $U$ has unique global minimum at $x^* \in M$. This is given by Proposition \ref{prop:convexity}, in Section \ref{sec:symm} below. 

\indent Proposition \ref{prop:convexity} motivates the main idea of the present work\,: the Riemannian barycentre $\bar{x}_{\scriptscriptstyle T}$ of $P_{\scriptscriptstyle{T}}$ can be used as a proxy for the global minimum $x^*$ of $U$. In~general, $\bar{x}_{\scriptscriptstyle T}$ only provides an approximation of $x^*$, but the two are equal if $U$ is invariant by geodesic symmetry about $x^*$, as stated in Proposition \ref{prop:main}, in Section \ref{sec:algorithm} below.  

The following Section \ref{sec:concentration} introduces Proposition \ref{sec:concentration}, which estimates the Riemannian distance between $\bar{x}_{\scriptscriptstyle T}$ and $x^*\,$, as a function of $T$. 

\section{Concentration of the barycentre} \label{sec:concentration}
Let $P$ be a probability distribution on a complete Riemannian manifold $M$. A (Riemannian) barycentre of $P$ is any global minimiser $\bar{x} \in M$ of the function
\begin{equation} \label{eq:variance}
\hspace{-0.5cm}  \mathcal{E}(x) = \frac{1}{2}\, \int_M d^{\scriptscriptstyle 2}(x,z)P(dz) \hspace{0.25cm}\mbox{ for }\, x \in M
\end{equation}
The following statement is due to Karcher, and was improved upon by Afsari~\cite{karcher}\cite{afsari}\,: \textit{if $P$ is supported inside a geodesic ball $B(x^*,\delta)$, where $x^* \in M$ and $\delta < \frac{1}{2} r_{\scriptscriptstyle cx}$ ($r_{\scriptscriptstyle cx}$ the convexity radius of $M$), then $\mathcal{E}$ is strongly convex on $B(x^*,\delta)$, and $P$ has a unique barycentre $\bar{x} \in B(x^*,\delta)$.} 

On the other hand, the present work considers a setting where $P$ is not supported inside $B(x^*,\delta)$, but merely concentrated on this ball. Precisely, assume $P$ is equal to the Gibbs distribution
\begin{equation} \label{eq:gibbs}
\hspace{-0.3cm} P_{\scriptscriptstyle{T}}(dz) \,=\, \left(Z(T)\right)^{-1}\,\exp\left[-\frac{U(z)}{T}\right]\mathrm{vol}(dz) \,\,;\, T > 0
\end{equation}
where $Z(T)$ is a normalising constant, $U$ is a $C^2$ function with unique global minimum at $x^*$, and $\mathrm{vol}$ is the Riemannian volume of $M$. Then, let $\mathcal{E}_{\scriptscriptstyle{T}}$ denote the function $\mathcal{E}$ in (\ref{eq:variance}), and let $\bar{x}_{\scriptscriptstyle{T}}$ denote any barycentre of $P_{\scriptscriptstyle{T\,}}$. 

In this new setting, it is not clear whether $\mathcal{E}_{\scriptscriptstyle{T}}$ is differentiable or not. Therefore, statements about convexity of $\mathcal{E}_{\scriptscriptstyle{T}}$ and uniqueness of $\bar{x}_{\scriptscriptstyle{T}}$ are postponed to the following Section \ref{sec:symm}. For now, it is possible to state the following Proposition \ref{prop:concentration}. In this proposition, $d(\cdot,\cdot)$ denotes Riemannian distance, and $W(\cdot,\cdot)$ denotes the Kantorovich ($L^1$-Wasserstein) distance~\cite{kantorovich}\cite{villani}. Moreover, $(\,\mu_{\scriptscriptstyle \min\,},\mu_{\scriptscriptstyle \max})$ is any open interval which contains the spectrum of the Hessian $\nabla^2 U(x^*)$, considered as a linear mapping of the tangent space $T_{\scriptscriptstyle x^*}M$.  

\begin{proposition} \label{prop:concentration}
 assume $M$ is an $n$-dimensional compact Riemannian manifold with non-negative sectional curvature. Denote $\delta_{\scriptscriptstyle x^*}$ the Dirac distribution at $x^*$. The following hold, \\[0.1cm]
\textbf{(i)} for any $\eta > 0$, 
\begin{equation} \label{eq:concentration1}
W(P_{\scriptscriptstyle{T\,}},\delta_{\scriptscriptstyle x^*}) < \frac{\eta^2}{(4\,\mathrm{diam}\, M)}\,\,\Longrightarrow\,\,
d(\bar{x}_{\scriptscriptstyle{T\,}},x^*) < \eta
\end{equation}
\textbf{(ii)} for $T \leq T_o$ (which can be computed explicitly)
\begin{equation} \label{eq:concentration2}
  W(P_{\scriptscriptstyle{T\,}},\delta_{\scriptscriptstyle x^*}) \leq \sqrt{2\pi}\,\left(\pi/2\right)^{n-1}\,B^{-1}_n\,\left(\mu_{\scriptscriptstyle \max}/\mu_{\scriptscriptstyle \min}\right)^{n/2}\,\left(T/\mu_{\scriptscriptstyle \min}\right)^{1/2}
\end{equation}
where $B_n = B(1/2,n/2)$ in terms of the Beta function. 
\end{proposition}
Proposition \ref{prop:concentration} is motivated by the idea of using $\bar{x}_{\scriptscriptstyle T}$ as an approximation of $x^*$. Intuitively, this requires choosing $T$ so small that $P_{\scriptscriptstyle T}$ is sufficiently close to $\delta_{\scriptscriptstyle x^*\,}$. Just how small a $T$ may be required is indicated by the inequality in (\ref{eq:concentration2}). This inequality is optimal and explicit, in the following sense.

It is optimal because the dependence on $T^{1/2}$ in its right-hand side cannot be improved. Indeed, by the multi-dimensional Laplace approximation (see~\cite{laplace}, for example), the left-hand side is equivalent  to $\mathrm{L}\cdot T^{1/2\,}$ (in the limit $T \rightarrow 0$). While this constant $\mathrm{L}$ is not tractable, the constants appearing in Inequality (\ref{eq:concentration2}) depend explicitly on the manifold $M$ and the function $U$. In fact, this inequality does not follows from the multi-dimensional Laplace approximation, but rather from volume comparison theorems of Riemannian geometry~\cite{chavel}. 

In spite of these nice properties, Inequality (\ref{eq:concentration2}) does not escape the curse of dimensionality. Indeed, for fixed $T$, its right-hand side increases exponentially with the dimension $n$ (note that $B_n$ decreases like $n^{\scriptscriptstyle -1/2}$). On the other hand, although $T_o$ also depends on $n$, it is typically much less affected by dimensionality, and decreases slower that $n^{-1}$ as $n$ increases.   

\section{Convexity and uniqueness} \label{sec:symm}
Assume now that $M$ is a simply-connected, compact Riemannian symmetric space. In this case, for any $T$, the function $\mathcal{E}_{\scriptscriptstyle T}$ turns out to be $C^2$ throughout $M$. This results from the following lemma. 
\begin{lemma} \label{eq:lemma}
 let $M$ be a simply-connected compact Riemannian symmetric space. Let $\gamma : I \rightarrow M$ be a geodesic defined on a compact interval $I$. Denote $\mathrm{Cut}(\gamma)$ the union of all cut loci $\mathrm{Cut}(\gamma(t))$ for $t \in I$.  Then, the topological dimension of $\mathrm{Cut}(\gamma)$ is strictly less than $n = \dim M$. In particular, $\mathrm{Cut}(\gamma)$ is a set with volume equal to zero. 

\textbf{Remark\,:} the assumption that $M$ is simply-connected cannot be removed, as the conclusion does not hold if $M$ is a real projective space.
\end{lemma}

The proof of Lemma \ref{eq:lemma} uses the structure of Riemannian symmetric spaces, as well as some results from topological dimension theory~\cite{helgason} (Chapter VII). The notion of topological dimension arises because it is possible $\mathrm{Cut}(\gamma)$ is not a manifold. The lemma immediately implies, for all $t$,
$$
\mathcal{E}_{\scriptscriptstyle T}(\gamma(t)) \,=\, \frac{1}{2}\,\int_{M}\,d^{\scriptscriptstyle 2}(\gamma(t),z)P_{\scriptscriptstyle T}(dz)\,=\, \frac{1}{2}\,\int_{M - \mathrm{Cut}(\gamma)}\,d^{\scriptscriptstyle 2}(\gamma(t),z)P_{\scriptscriptstyle T}(dz)
$$
Then, since the domain of integration avoids the cut loci of all the $\gamma(t)$, it becomes possible to differentiate under the integral. This is used in obtaining the following (the assumptions are the same as in Lemma \ref{eq:lemma}).
\begin{corollary}\label{corr}
for $x \in M$, let $G_x(z)=\nabla f_z(x)$ and $H_x(z) = \nabla^2f_z(x)$, where $f_z$ is the function $x\mapsto \frac{1}{2}\,d^{\scriptscriptstyle2}(x,z)$. The following integrals converge for any $T$
$$
G_x \,=\, \int_{M-\mathrm{Cut}(x)}\,G_x(z)\,P_{\scriptscriptstyle T}(dz) \hspace{0.25cm};\hspace{0.25cm}
H_x \,=\, \int_{M-\mathrm{Cut}(x)}\,H_x(z)\,P_{\scriptscriptstyle T}(dz)
$$
and both depend continuously on $x$. Moreover,  
\begin{equation} \label{eq:derivatives}
\nabla \mathcal{E}_{\scriptscriptstyle T}(x) = G_x \hspace{0.2cm} \mbox{and}\hspace{0.2cm}
\nabla^2 \mathcal{E}_{\scriptscriptstyle T}(x) = H_x
\end{equation}
so that $\mathcal{E}_{\scriptscriptstyle T}$ is $C^2$ throughout $M$.
\end{corollary}
With Corollary  \ref{corr} at hand, it is possible to obtain Proposition \ref{prop:convexity}, which is concerned with the convexity of $\mathcal{E}_{\scriptscriptstyle T}$ and uniqueness of $\bar{x}_{\scriptscriptstyle T\,}$. In this proposition, the following notation is used
\begin{equation} \label{eq:fnt}
   f(T) =  (2/\pi)\left(\pi/8\right)^{n/2}\left(\mu_{\scriptscriptstyle \max}/T\right)^{n/2}\exp\left(-U_{\scriptscriptstyle\delta}/T\right)
\end{equation}
where $U_{\scriptscriptstyle \delta} = \inf \lbrace U(x) - U(x^*)\,;\, x \notin B(x^*,\delta)\rbrace$ for positive $\delta$. The reader may wish to note the fact that $f(T)$ decreases to $0$ as $T$ decreases to $0$. 
\begin{proposition} \label{prop:convexity}
let $M$ be a simply-connected compact Riemannian symmetric space. Let $\kappa^2$ be the maximum sectional curvature of $M$, and $r_{\scriptscriptstyle cx} = \kappa^{-1}\frac{\pi}{2}$ its convexity radius. If $T \leq T_o$ (see \textit{(ii)} of Proposition \ref{prop:concentration}), then the following hold for any $\delta < \frac{1}{2} r_{\scriptscriptstyle cx}$. 
\\[0.1cm]
\textbf{(i)} for all $x$ in the geodesic ball $B(x^*,\delta)$, 
\begin{equation} \label{eq:convexity1}
\nabla^2 \mathcal{E}_{\scriptscriptstyle T}(x) \geq  \mathrm{Ct}(2\delta)\left( 1 - \mathrm{vol}(M) f(T)\right) - \pi  A_M f(T)
\end{equation}
where $\mathrm{Ct}(2\delta) = 2\kappa\delta\cot(2\kappa\delta) > 0$ and $A_M > 0$ is a constant given by the structure of the symmetric space $M$. \\[0.1cm]
\textbf{(ii}) there exists $T_{\scriptscriptstyle \delta}$ (which can be computed explicitly), such that $T \leq T_{\scriptscriptstyle \delta}$ implies $\mathcal{E}_{\scriptscriptstyle{T}}$ is strongly convex on $B(x^*,\delta)\,$, and has a unique global minimum $\bar{x}_{\scriptscriptstyle T} \in B(x^*,\delta)$. In particular, this means $\bar{x}_{\scriptscriptstyle T}$ is the unique barycentre of $P_{\scriptscriptstyle T\,}$. 
\end{proposition}
Note that \textit{(ii)} of Proposition \ref{prop:convexity} generalises the statement due to Karcher~\cite{karcher}, which was recalled in Section \ref{sec:concentration}. 
\section{Finding $T_o$ and $T_{\scriptscriptstyle \delta}$} \label{sec:ts}
Propositions \ref{prop:concentration} and \ref{prop:convexity} claim that $T_o$ and $T_{\scriptscriptstyle \delta}$ can be computed explicitly.  This means 
that, with some knowledge of the Riemannian manifold $M$ and the function $U$, $T_o$ and $T_{\scriptscriptstyle \delta}$ can be found by solving scalar equations. The current section gives the definitions of $T_o$ and $T_{\scriptscriptstyle \delta\,}$.

In the notation of Proposition \ref{prop:concentration}, let $\rho > 0$ be small enough, so that, 
$$
\mu_{\scriptscriptstyle \min\,}d^{\scriptscriptstyle 2}(x,x^*) \,\leq\, 2\left(U(x)-U(x^*)\right)\,\leq\, \mu_{\scriptscriptstyle \max\,}d^{\scriptscriptstyle 2}(x,x^*)
$$
whenever $d(x,x^*) \leq \rho\,$, and consider the quantity
$$
f(T,m,\rho) \,=\, (2/\pi)^{1/2}\,\left(\mu_{\scriptscriptstyle\max}/T\right)^{m/2}\,\exp\left(-U_{\scriptscriptstyle\rho}/T\right)
$$
where $U_{\scriptscriptstyle \rho}$ is defined as in (\ref{eq:fnt}). Note that $f(T,m,\rho)$ decreases to $0$ as $T$ decreases to $0$, for fixed $m$ and $\rho$. Now, it is possible to define $T_o$ as
\begin{equation}
  T_o \,=\, \min\left\lbrace T^1_o\,,\,T^2_o\right\rbrace \hspace{0.25cm}\mbox{ where }
\end{equation}
$$
\begin{array}{ll}
T^1_o =& \inf\left\lbrace T> 0\,:\, f(T,n-2,\rho)\,>\,\rho^{2-n}\,A_{n-1}\,\right\rbrace \\[0.3cm]
T^2_o =& \inf\left\lbrace T> 0\,:\, f(T,n+1,\rho)\,>\,\left(\mu_{\scriptscriptstyle\max}/\mu_{\scriptscriptstyle\min}\right)^{n/2}\,C_n\right\rbrace \\[0.2cm]
\end{array}
$$
Here, $A_n = E|X|^n$ for $X \sim  N(0,1)$, and $C_n = \omega_n\,A_{n}/\!\left(\mathrm{diam}\, M\times \mathrm{vol}\,M\right)$, where $\omega_n$ is the surface area of a unit sphere $S^{n-1\,}$. \vfill\pagebreak

With regard to Proposition \ref{prop:convexity}, define $T_{\scriptscriptstyle \delta}$ as follows,
\begin{equation} \label{eq:camille}
  T_{\scriptscriptstyle \delta} \,=\, \min\left\lbrace T^1_{\scriptscriptstyle \delta}\,,\,T^2_{\scriptscriptstyle \delta}\right\rbrace - \varepsilon 
\end{equation}
for some arbitrary $\varepsilon > 0$. Here, in the notation of (\ref{eq:concentration2}), (\ref{eq:fnt}) and (\ref{eq:convexity1}), 
$$
\begin{array}{ll}
T^1_{\scriptscriptstyle \delta} =& \inf
\left\lbrace T\leq T_o\,:\, \sqrt{2\pi}\,(T/\mu_{\scriptscriptstyle\min})^{1/2}\,>\,\delta^2\,\left(\mu_{\scriptscriptstyle\min}/\mu_{\scriptscriptstyle\max}\right)^{n/2}\,D_n\right\rbrace \\[0.3cm]
T^2_{\scriptscriptstyle \delta} =& \inf\left\lbrace T\leq T_o \,:\, f(T)\,>\, \mathrm{Ct}(2\delta)\left(\mathrm{Ct}(2\delta)\,\mathrm{vol}(M) + \pi A_M\right)^{-1}\,\right\rbrace \\[0.2cm]
\end{array}
$$
where $D_n = \,(2/\pi)^{n-1}\,B_n/(4\,\mathrm{diam}\,M)$.
  
\section{Black-box optimisation} \label{sec:algorithm}
Consider the problem of searching for the unique global minimum $x^*$ of $U$. In~black-box optimisation, it is only possible to evaluate $U(x)$ for given
$x \in M$, and the cost of this evaluation precludes numerical approximation of derivatives. Then,~the problem is to find $x^*$ using successive evaluations of $U(x)$ (hopefully, as few of these evaluations as possible). 

Here, a new algorithm for solving this problem is described. The idea of this algorithm is to find $\bar{x}_{\scriptscriptstyle T}$ using successive evaluations of $U(x)$, in the hope that $\bar{x}_{\scriptscriptstyle T}$ will provide a good approximation of $x^*$. While the quality of this approximation is controlled by Inequalities (\ref{eq:concentration1}) and (\ref{eq:concentration2}) of Proposition \ref{prop:concentration}, in some cases of interest, $\bar{x}_{\scriptscriptstyle T}$ is exactly equal to $x^*$, for correctly chosen $T$, as in the following proposition~\ref{prop:main}.

To state this proposition, let $s_{\scriptscriptstyle{x^*}}$ denote geodesic symmetry about $x^*$ (see~\cite{helgason}). This is the transformation of $M$, which leaves $x^*$ fixed, and reverses the direction of geodesics passing through $x^*$. 
\begin{proposition} \label{prop:main}
 assume that $U$ is invariant by geodesic symmetry about $x^*\,$, in the sense that $U \circ s_{\scriptscriptstyle{x^*}} = U$. If $T \leq T_{\scriptscriptstyle \delta}$ (see (ii) of Proposition \ref{prop:convexity}), then $\bar{x}_{\scriptscriptstyle T} = x^*$ is the unique barycentre of $P_{\scriptscriptstyle T\,}$. 
\end{proposition}
Proposition \ref{prop:main} follows rather directly from Proposition \ref{prop:convexity}. Precisely, by \textit{(ii)} of Proposition~\ref{prop:convexity}, the condition $T \leq T_{\scriptscriptstyle \delta}$ implies $\mathcal{E}_{\scriptscriptstyle{T}}$ is strongly convex on $B(x^*,\delta)$, and $\bar{x}_{\scriptscriptstyle T} \in B(x^*,\delta)$. Thus, $\bar{x}_{\scriptscriptstyle T}$ is the unique stationary point of $\mathcal{E}_{\scriptscriptstyle{T}}$ in $B(x^*,\delta)$. But, using the fact that $U$ is invariant by geodesic symmetry about $x^*\,$, it is possible to prove that $x^*$ is a stationary point of $\mathcal{E}_{\scriptscriptstyle{T\,}}$, and this implies  $\bar{x}_{\scriptscriptstyle T} = x^*$.\\[0.1cm]
The two following examples verify the conditions of Proposition \ref{prop:main}.\\[0.1cm]
\indent \textbf{Example 1\,:} assume $M = \mathrm{Gr}(k,\mathbb{C}^n)$ is a complex Grassmann manifold. In particular, $M$ is a simply-connected, compact Riemannian symmetric space. Identify $M$ with the set of Hermitian projectors $x : \mathbb{C}^n\rightarrow \mathbb{C}^n$ such that $\mathrm{tr}(x) = k$, where $\mathrm{tr}$ denotes the trace. Then, define $U(x) = -\,\mathrm{tr}(C\,x)$ for $x \in \mathrm{Gr}(k,\mathbb{C}^n)$, where $C$ is a Hermitian positive-definite matrix with distinct eigenvalues. Now, the unique global minimum of $U$ occurs at $x^*$, the projector onto the principal \\ $k$-subspace of $C$. Also, the geodesic symmetry $s_{\scriptscriptstyle x^*}$ is given by $s_{\scriptscriptstyle x^*}\cdot x = r_{\scriptscriptstyle x^*}x\,r_{\scriptscriptstyle x^*\,}$, where $r_{\scriptscriptstyle x^*} : \mathbb{C}^n\rightarrow \mathbb{C}^n$ denotes reflection through the image space of $x^*$. It is elementary to verify that $U$ is invariant by this geodesic symmetry. \vfill\pagebreak
\indent \textbf{Example 2\,:} let $M$ be a simply-connected, compact Riemannian symmetric space, and $U_{\scriptscriptstyle o}$ a function on $M$ with unique global minimum at $o \in M$. Assume moreover that $U_{\scriptscriptstyle o}$ is invariant by geodesic symmetry about $o$. For each $x^* \in M$, there exists an isometry $g$ of $M$, such that $x^* = g\cdot o$. Then, $U(x) = U_{\scriptscriptstyle o}(g^{\scriptscriptstyle{-1}}\cdot x)$ has unique global minimum at $x^*$, and is invariant by geodesic symmetry about~$x^*$. \\[0.1cm]
\indent Example 1 describes the standard problem of finding the principal subspace of the covariance matrix $C$. In Example 2, the function $U_{\scriptscriptstyle o}$ is a known template, which undergoes an unknown transformation $g$, leading to the observed pattern~$U$. This is a typical situation in pattern recognition problems. 

Of course, from a mathematical point of view, Example 2 is not really an example, since it describes the completely general setting where the conditions of Proposition \ref{prop:main} are verified. In this setting, consider the following algorithm. \\[0.2cm]
\indent \textbf{Description of the algorithm\,:} \\[0.1cm]
-- input\,: $T \leq T_{\scriptscriptstyle \delta}$ \hfill \% to find such $T$, see Section \ref{sec:ts}\hspace{0.524cm}   \\[0.05cm]
$\phantom{\mbox{-- input\,:}}$ $Q(x,dz) = q(x,z)\mathrm{vol}(dz)$         \hfill \% symmetric Markov kernel $\phantom{\mbox{on $M$}}$   \\[0.05cm]
$\phantom{\mbox{-- input\,:}}$ $\hat{x}_{\scriptscriptstyle 0} = z_{\scriptscriptstyle 0} \in M$       \hfill \% initial guess for $x^*$ \hspace{1.93cm}  \\[0.1cm]
-- iterate\,: for $n = 1, 2, \ldots$ \\[0.05cm]
$\phantom{\mbox{-- iterate\,:}}$ \texttt{(1)} sample $z_{\scriptscriptstyle n} \sim q(z_{\scriptscriptstyle n-1},z)$  \\[0.15cm]
$\phantom{\mbox{-- iterate\,:}}$ \texttt{(2)} compute $r_{\scriptscriptstyle n} = 1 - \min\left\lbrace 1,\exp\left[\left(U(z_{\scriptscriptstyle n-1}) - U(z_{\scriptscriptstyle n})\right)\middle/ T\right]\right\rbrace$ \\[0.15cm]
$\phantom{\mbox{-- iterate\,:}}$ \texttt{(3)} reject $z_{\scriptscriptstyle n}$ with probability $r_{\scriptscriptstyle n}$ \hfill \% then, $z_{\scriptscriptstyle n} = 
z_{\scriptscriptstyle n-1}$ \hspace{2.42cm} \\[0.1cm]
$\phantom{\mbox{-- iterate\,:}}$ \texttt{(4)} $\hat{x}_{\scriptscriptstyle n} = \hat{x}_{\scriptscriptstyle n-1}\, \#_{\scriptscriptstyle \frac{1}{n}}\, z_{\scriptscriptstyle n}$ \hfill \% see definition (\ref{eq:weighted}) below \hspace{0.99cm} \\[0.15cm]
-- until\,: $\hat{x}_{\scriptscriptstyle n}$ does not change sensibly \\[0.1cm]
-- output\,: $\hat{x}_{\scriptscriptstyle n}$ \hfill \% approximation of $x^*$ \hspace{1.66cm}  \\[0.3cm]
The above algorithm recursively computes the Riemannian barycentre $\hat{x}_{\scriptscriptstyle n\,}$ of the samples $z_{\scriptscriptstyle n}$ generated by a symmetric Metropolis-Hastings algorithm (see~\cite{rosenthal}). Here, The Metropolis-Hastings algorithm is implemented in lines \texttt{(1)--(3)}. On~the other hand, line \texttt{(4)} takes care of the Riemannian barycentre. Precisely,
if $\gamma:[0,1]\rightarrow M$ is a length-minimising geodesic connecting $\hat{x}_{\scriptscriptstyle n-1}$ to $z_{\scriptscriptstyle n\,}$, let 
\begin{equation} \label{eq:weighted}
  \hat{x}_{\scriptscriptstyle n-1}\, \#_{\scriptscriptstyle \frac{1}{n}}\, z_{\scriptscriptstyle n} \,=\,\gamma\left(1/n\right)
\end{equation}
This geodesic $\gamma$ need not be unique. 

The point of using the Metropolis-Hastings algorithm is that the generated $z_{\scriptscriptstyle n}$ eventually sample from the Gibbs distribution $P_{\scriptscriptstyle T\,}$. The convergence of the distribution $P_{\scriptscriptstyle n}$ of $z_{\scriptscriptstyle n}$ to $P_{\scriptscriptstyle T}$ takes place exponentially fast. Indeed, 
it may be inferred from~\cite{rosenthal} (see Theorem 8, Page 36)
\begin{equation} \label{eq:ergo}
\Vert P_{\scriptscriptstyle n} - P_{\scriptscriptstyle T}\Vert_{\scriptscriptstyle TV} \leq (1-p_{\,\scriptscriptstyle T}\,)^n 
\end{equation}
where $\Vert \cdot\Vert_{\scriptscriptstyle TV}$ is the total variation norm, and $p_{\,\scriptscriptstyle T} \in (0,1)$ verifies
$$
p_{\,\scriptscriptstyle T} \leq \, (\mathrm{vol}\,M)\,\inf_{x,z} \,q(x,z)\,\exp(-\sup_x U(x)/T)\, 
$$
so the rate of convergence is degraded when $T$ is small. 
 
Accordingly, the intuitive justification of the above algorithm is the following. Since the $z_{\scriptscriptstyle n}$ eventually sample from the Gibbs distribution $P_{\scriptscriptstyle T\,}$, and the desired global minimum $x^*$ of $U$ is equal to the barycentre $\bar{x}_{\scriptscriptstyle T}$ of $P_{\scriptscriptstyle T\,}$ (by Proposition \ref{prop:main}), then the barycentre $\hat{x}_{\scriptscriptstyle n}$ of the $z_{\scriptscriptstyle n}$ is expected to converge to $x^*$. 

It should be emphasised that, in the present state of the literature, there is no rigorous result which confirms this convergence $z_{\scriptscriptstyle n} \rightarrow x^*\,$. It is therefore an open problem, to be confronted in future work.  

For a basic computer experiment, consider $M = S^2 \subset \mathbb{R}^3,$ and let
\begin{equation} \label{eq:example1}
  U(x) = -\,P_{\scriptscriptstyle 9}(x^{\scriptscriptstyle 3}) \hspace{0.25cm}\mbox{ for }\, x = (x^{\scriptscriptstyle 1},x^{\scriptscriptstyle 2},x^{\scriptscriptstyle 3}) \in S^2
\end{equation}
where $P_{\scriptscriptstyle 9}$ is the Legendre polynomial of degree $9$~\cite{wongbeals}. The unique global minimiser of $U$ is $x^* = (0,0,1)$, and the conditions of Proposition \ref{prop:main} are verified, since $U$ is invariant by reflection in the $x^{\scriptscriptstyle 3}$ axis, which is geodesic symmetry about $x^*$. 
\vspace{-3.1cm}
\begin{figure}
  \begin{minipage}[h]{0.5\textwidth}
    \includegraphics[width=\textwidth]{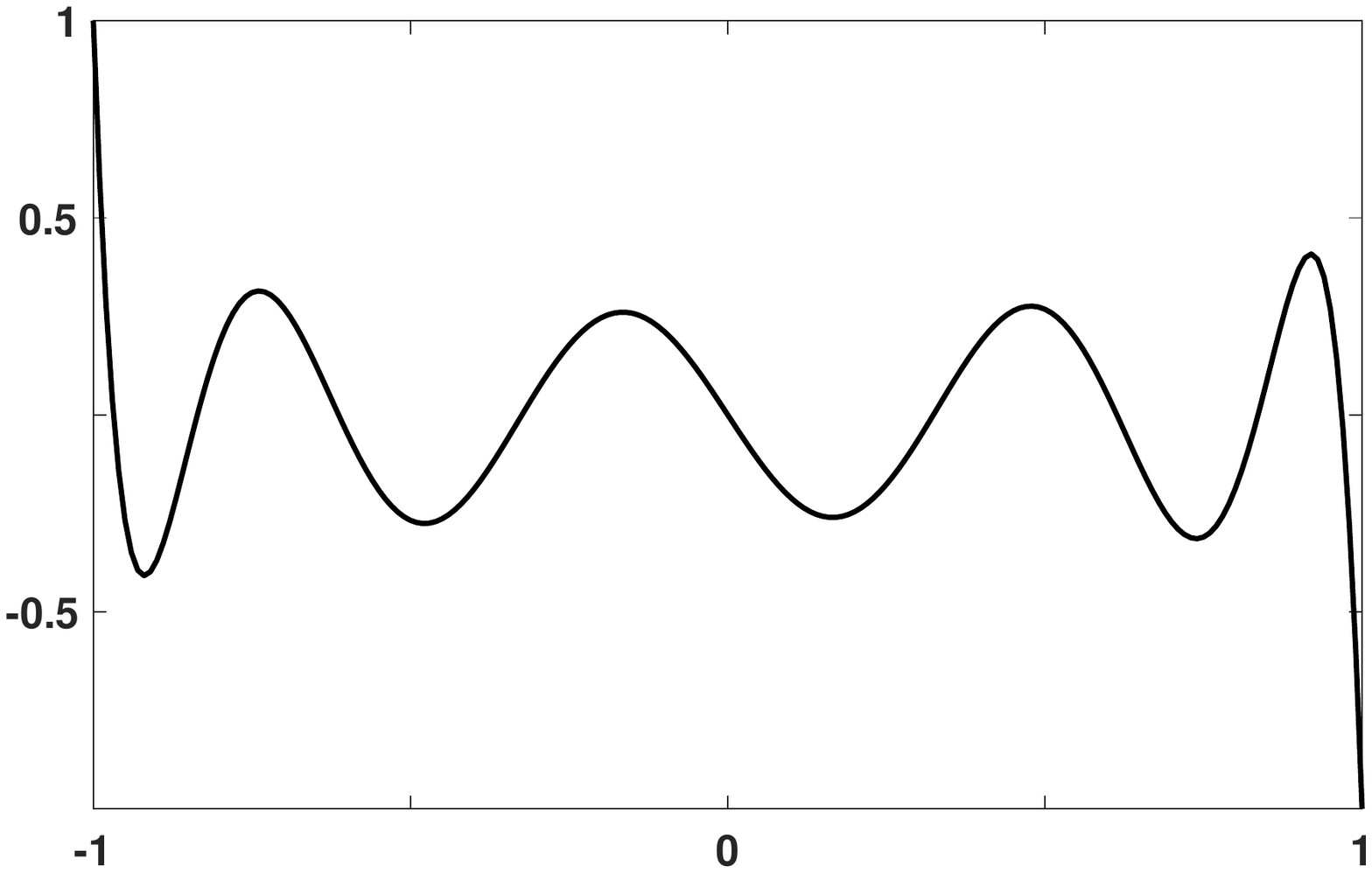}
\vspace{-3.15cm}
    \caption{graph of $-P_{\scriptscriptstyle 9}(x^{\scriptscriptstyle 3})$}
    \label{fig:1}
  \end{minipage}
\hspace{0.1cm}
  \begin{minipage}[h]{0.5\textwidth}
    \includegraphics[width=\textwidth]{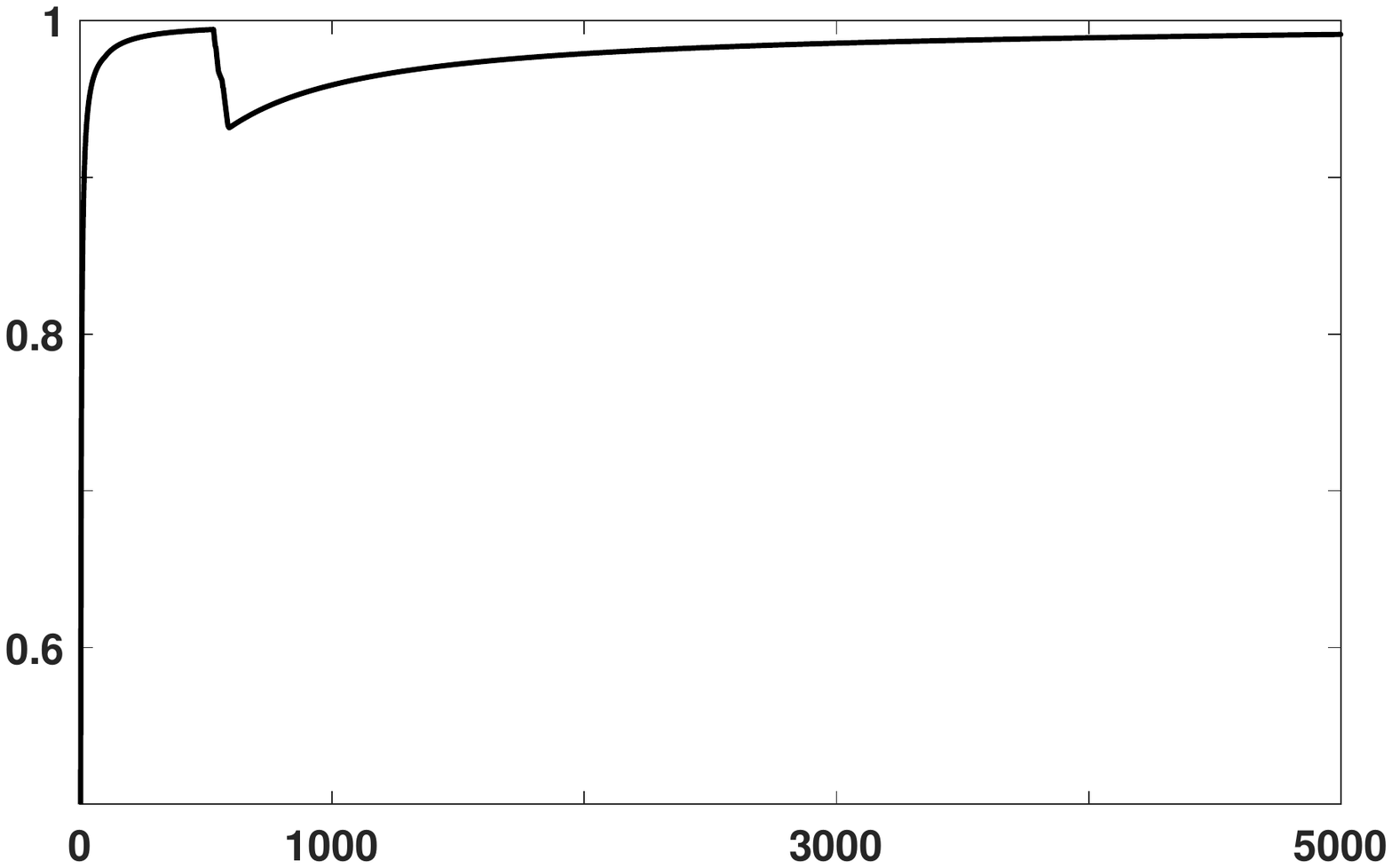}
\vspace{-3.15cm}
    \caption{$\hat{x}^{\scriptscriptstyle 3}_{\scriptscriptstyle n}$ versus $n$}
   \label{fig:2}
  \end{minipage}
\end{figure}

\vspace{-0.33cm}
Figure \ref{fig:1} shows the dependence of $U(x)$ on $x^{\scriptscriptstyle 3}$, displaying multiple local minima and maxima. Figure \ref{fig:2} shows the algorithm overcoming these local minima and maxima, and converging to the global minimum $x^* = (0,0,1)$, within $n=5000$ iterations. The experiment was conducted with $T = 0.2$, and the Markov kernel $Q$ obtained from the von Mises-Fisher distribution (see~\cite{mardia}).  The initial guess 
$\hat{x}_{\scriptscriptstyle 0} = (0,0,-1)$ is not shown in Figure \ref{fig:2}. 

In comparison, a standard simulated annealing method offered less robust performance, which varied considerably with the choice of annealing schedule. 

\vfill
\pagebreak

\section{Proofs}
This section is devoted to the proofs of the results stated in previous sections. 

As of now, assume that $U(x^*) = 0$. There is nos loss of generality in making this assumption. 

\subsection{Proof of Proposition \ref{prop:concentration}}
\begin{subequations}
\textbf{Proof of \textit{(i)}\,:} denote $f_x(z) = \frac{1}{2}\,d^{\scriptscriptstyle2}(x,z)$\,. By the definition of $\mathcal{E}_{\scriptscriptstyle{T}}$
\begin{equation} \label{eq:prcon1}
  \mathcal{E}_{\scriptscriptstyle{T}}(x) \,=\, \int_M\, f_x(z)\,P_{\scriptscriptstyle{T}}(dz)
\end{equation}
Moreover, let $\mathcal{E}_{\scriptscriptstyle{0}}$ be the function
\begin{equation} \label{eq:prcon2}
 \mathcal{E}_{\scriptscriptstyle{0}}(x) \,=\, \int_M\, f_x(z)\,\delta_{\scriptscriptstyle{x^*}}(dz) \,=\, \frac{1}{2}\,d^{\scriptscriptstyle 2}(x,x^*)
\end{equation}
For any $x$, it is elementary that $f_x(z)$ is Lipschitz continuous, with respect to $z$, with Lipschitz constant $\mathrm{diam}\, M$. Then, from the Kantorovich-Rubinshtein formula~\cite{villani}, 
\begin{equation} \label{eq:prcon3}
  \left\vert \mathcal{E}_{\scriptscriptstyle{T}}(x)\,-\, \mathcal{E}_{\scriptscriptstyle{0}}(x)\right\vert \leq\,\left(\mathrm{diam}\, M\right)W(P_{\scriptscriptstyle{T\,}},\delta_{\scriptscriptstyle x^*})
\end{equation}
a uniform bound in $x \in M$. It now follows that
\begin{equation} \label{eq:prcon4}
\inf_{x\in B(x^*,\eta)}\, \mathcal{E}_{\scriptscriptstyle{T}}(x) \,\,- \inf_{x\in B(x^*,\eta)}\, \mathcal{E}_{\scriptscriptstyle{0}}(x)\,\,\leq \,\left(\mathrm{diam}\, M\right)W(P_{\scriptscriptstyle{T\,}},\delta_{\scriptscriptstyle x^*}) \hspace{0.2cm}\mbox{ and}
\end{equation}
\begin{equation} \label{eq:prcon5}
\inf_{x\notin B(x^*,\eta)}\, \mathcal{E}_{\scriptscriptstyle{0}}(x) \,\,- \inf_{x\notin B(x^*,\eta)}\, \mathcal{E}_{\scriptscriptstyle{T}}(x)\,\,\leq \,\left(\mathrm{diam}\, M\right)W(P_{\scriptscriptstyle{T\,}},\delta_{\scriptscriptstyle x^*}) \hspace{0.2cm}\phantom{\mbox{ and}}
\end{equation}
However, from (\ref{eq:prcon2}), it is clear that
$$
\inf_{x\in B(x^*,\eta)}\, \mathcal{E}_{\scriptscriptstyle{0}}(x) = 0 \;\;\;\;\;\mbox{and}\;\,\,\,
\inf_{x\notin B(x^*,\eta)}\, \mathcal{E}_{\scriptscriptstyle{0}}(x) = \frac{\eta^2}{2}
$$
To complete the proof, replace this into (\ref{eq:prcon4}) and (\ref{eq:prcon5}). Then,  assuming the condition in (\ref{eq:concentration1}) is verified,
\begin{equation}
 \inf_{x\in B(x^*,\eta)}\, \mathcal{E}_{\scriptscriptstyle{T}}(x) \,< \frac{\eta^2}{4} <\,  \inf_{x\notin B(x^*,\eta)}\, \mathcal{E}_{\scriptscriptstyle{T}}(x)
\end{equation}
This means that any global minimum $\bar{x}_{\scriptscriptstyle T}$ of $\mathcal{E}_{\scriptscriptstyle{T}}$ must belong to the open ball $B(x^*,\eta)$. In other words, $d(\bar{x}_{\scriptscriptstyle T\,},x^*) < \eta$. This completes the proof of (\ref{eq:concentration1}). \hfill$\blacksquare$ \\[0.2cm]
\end{subequations}
\textbf{Proof of \textit{(ii)}\,:} let $\rho \leq \min\lbrace \mathrm{inj}\,x^*,\,\kappa^{-1}\,\frac{\pi}{2}\rbrace$ where $\mathrm{inj}\,x^*$ is the injectivity radius of $M$ at $x^*$, and $\kappa^2$ is an upper bound on the sectional curvature of $M$. Assume, in addition, $\rho$ is small enough so 
\begin{subequations}
\begin{equation} \label{eq:prcon21}
\mu_{\scriptscriptstyle \min\,}d^{\scriptscriptstyle 2}(x,x^*) \,\leq\, 2\left(U(x)-U(x^*)\right)\,\leq\, \mu_{\scriptscriptstyle \max\,}d^{\scriptscriptstyle 2}(x,x^*)
\end{equation}
whenever $d(x,x^*) \leq \rho\,$. Further, consider the truncated distribution
\begin{equation} \label{eq:prcon22}
  P^\rho_{\scriptscriptstyle T}(dz) \,=\, \frac{\mathbf{1}_{\scriptscriptstyle {B_\rho}}(z)}{P_{\scriptscriptstyle T}(B_\rho)}\cdot P_{\scriptscriptstyle T}(dz)
\end{equation}
where $\mathbf{1}$ denotes the indicator function, and $B_\rho$ stands for the open ball $B(x^*,\rho)$. 
Of course, by the triangle inequality,
\begin{equation} \label{eq:prcon23}
  W(P_{\scriptscriptstyle T},\delta_{\scriptscriptstyle x^*}) \leq 
  W(P_{\scriptscriptstyle T},P^\rho_{\scriptscriptstyle T}) +
  W(P^\rho_{\scriptscriptstyle T},\delta_{\scriptscriptstyle x^*}) 
\end{equation} 
The proof relies on the following estimates, which use the notation of Section \ref{sec:ts}.\\[0.1cm]
\indent \textbf{First estimate\,:} if $T \leq T^1_o\,$, then
\begin{equation} \label{eq:estimate1}
W(P_{\scriptscriptstyle T},P^\rho_{\scriptscriptstyle T}) \leq \left(\mathrm{diam}\, M\times \mathrm{vol}\,M\right)\,\frac{2}{\pi}\,
\left(\frac{\pi}{8}\right)^{n/2}\left( \frac{\mu_{\scriptscriptstyle \max}}{T}\right)^{n/2}\exp\,\left(-U_{\scriptscriptstyle\rho}/T\right)
\end{equation}
\indent \textbf{Second estimate\,:} if $T \leq T^1_o\,$, then
\begin{equation} \label{eq:estimate2}
W(P^\rho_{\scriptscriptstyle T},\delta_{\scriptscriptstyle x^*}) \leq \,
2\,\sqrt{2\pi}\,\left(\frac{\pi}{2}\right)^{n-1}\,B^{-1}_n\,\left(\frac{\mu_{\scriptscriptstyle \max}}{\mu_{\scriptscriptstyle \min}}\right)^{n/2}\,\left(\frac{T}{\mu_{\scriptscriptstyle \min}}\right)^{1/2}
\end{equation}
These two estimates are proved below. Assume now they hold true, and $T \leq T_o\,$. In particular, since $T \leq T^2_o\,$, the definition of $T^2_o$ implies
$$
f(T,n+1,\rho) \leq \left(\mu_{\scriptscriptstyle\max}/\mu_{\scriptscriptstyle\min}\right)^{n/2}\,C_n
$$
Recall the definition of $C_n\,$, and express $\omega_n$ and $A_n$ in terms of the Gamma function~\cite{wongbeals}. The last inequality becomes
$$
\left(\mathrm{diam}\, M\times \mathrm{vol}\,M\right)\,f(T,n+1,\rho) \leq 2\,(2\pi)^{n/2}B^{-1}_n\,\left(\mu_{\scriptscriptstyle\max}/\mu_{\scriptscriptstyle\min}\right)^{n/2}
$$
This is the same as
$$
\left(\mathrm{diam}\, M\times \mathrm{vol}\,M\right)\,\frac{1}{\pi}\,\left(\frac{\pi}{8}\right)^{n/2}f(T,n+1,\rho) \leq \left(\frac{\pi}{2}\right)^{n-1}\,B^{-1}_n\,\left(\mu_{\scriptscriptstyle\max}/\mu_{\scriptscriptstyle\min}\right)^{n/2}
$$
By the definition of $f(T,n+1,\rho)$, it now follows the right-hand side of (\ref{eq:estimate1}) is less than half the right-hand side of (\ref{eq:estimate2}). 
In this case, (\ref{eq:concentration2}) follows from the triangle inequality (\ref{eq:prcon23}). \hfill$\blacksquare$ \\[0.2cm]
\end{subequations}
\begin{subequations}  
\indent \textbf{Proof of first estimate\,:} consider the coupling of $P_{\scriptscriptstyle T}$ and $P^\rho_{\scriptscriptstyle T\,}$, provided by the probability distribution $K$ on $M \times M$, 
\begin{equation} \label{eq:coupling}
K(dz_{\scriptscriptstyle 1}\times dz_{\scriptscriptstyle 2}) = P^\rho_{\scriptscriptstyle T}(dz_{\scriptscriptstyle 1})
\left[P_{\scriptscriptstyle T}(B_\rho)\delta_{ z_{\scriptscriptstyle 1}}(dz_{\scriptscriptstyle 2})\,+\,\mathbf{1}_{\scriptscriptstyle {B^c_\rho}}(z_{\scriptscriptstyle2})P_{\scriptscriptstyle T}(dz_{\scriptscriptstyle 2}) \right]
\end{equation}
where $B^c_\rho$ denotes the complement of $B_\rho\,$. Recall the definition of the Kantorovich distance (see~\cite{villani}). Replacing (\ref{eq:coupling}) into this definition, it follows that
\begin{equation} \label{eq:prest11}
  W(P_{\scriptscriptstyle T},P^\rho_{\scriptscriptstyle T}) \leq \left(\mathrm{diam}\,M\right)\,P_{\scriptscriptstyle T}(B^c_\rho)
\end{equation}
Then, from the definition (\ref{eq:gibbs}) of $P_{\scriptscriptstyle T\,}$,
\begin{equation} \label{eq:prest12}
P_{\scriptscriptstyle T}(B^c_\rho) \leq \left(Z(T)\right)^{-1}\left(\mathrm{vol}\,M\right)\,\exp\left(-U_{\scriptscriptstyle\rho}/T\right)
\end{equation}
Now, (\ref{eq:estimate1}) follows directly from (\ref{eq:prest11}) and (\ref{eq:prest12}), if the following lower bound on $Z(T)$ can be proved,
\begin{equation} \label{eq:lower}
Z(T) \geq \frac{\pi}{2}\,\left(\frac{8}{\pi}\right)^{n/2}\left( \frac{T}{\mu_{\scriptscriptstyle \max}}\right)^{n/2} \hspace{0.25cm}\mbox{ for } T\leq T^1_o
\end{equation}
To prove this lower bound, note that 
$$
Z(T) \,=\, \int_M\,e^{-\frac{U(z)}{\mathstrut T}}\,\mathrm{vol}(dz) \,\,\geq\, 
\int_{B_\rho}\,e^{-\frac{U(z)}{\mathstrut T}}\,\mathrm{vol}(dz)
$$
Using this last inequality and (\ref{eq:prcon21}), it is possible to write
\begin{equation} \label{eq:later}
Z(T) \geq\, 
\int_{B_\rho}\,e^{-\frac{U(z)}{\mathstrut T}}\,\mathrm{vol}(dz)\geq\, 
\int_{B_\rho}\,e^{-\frac{\mu_{\scriptscriptstyle \max}}{\mathstrut 2T}\,d^{\scriptscriptstyle 2}(x,x^*)}\,\,\mathrm{vol}(dz)
\end{equation}
Writing this last integral in Riemannian spherical coordinates, 
\begin{equation} \label{eq:comparison11}
  \int_{B_\rho}\,e^{-\frac{\mu_{\scriptscriptstyle \max}}{\mathstrut 2T}\,d^{\scriptscriptstyle 2}(x,x^*)}\,\,\mathrm{vol}(dz) \,=\,
  \int^\rho_0\int_{S^{n-1}}\,e^{-\frac{\mu_{\scriptscriptstyle \max}}{\mathstrut 2T}\,r^{\scriptscriptstyle 2}}\,\lambda(r,s)\,dr\,\omega_n(ds)
\end{equation}
where $\lambda(r,s)$ is the volume density in the Riemannian spherical coordinates, $r\geq 0$ and $s \in S^{n-1}$, and where $\omega_n(ds)$ is the area element of $S^{n-1}$. From the volume comparison theorem in~\cite{chavel} (see Page 129), 
\begin{equation} \label{eq:comparison12}
\lambda(r,s) \geq \left(\kappa^{-1}\sin(\kappa\,r)\right)^{\!n-1} \geq \left( (2/\pi)\,r\right)^{n-1}
\end{equation}
where the second inequality follows since $x\mapsto \sin(x)$ is concave for $x \in (0,\pi)$. Now, it follows from (\ref{eq:later}) and (\ref{eq:comparison11}),
\begin{equation} \label{eq:comparison13}
Z(T) \geq \,\omega_n\,\left(\frac{2}{\pi}\right)^{\!n-1}\,
\int^\rho_0\,e^{-\frac{\mu_{\scriptscriptstyle \max}}{\mathstrut 2T}\,r^{\scriptscriptstyle 2}}\,r^{n-1}\,\,dr
\end{equation}
where $\omega_n$ is the surface area of $S^{n-1\,}$. Thus, the required lower bound (\ref{eq:lower}) follows by noting that
$$
\!\! \int^\rho_0\,e^{-\frac{\mu_{\scriptscriptstyle \max}}{\mathstrut 2T}\,r^{\scriptscriptstyle 2}}\,r^{n-1}\,\,dr \,=\,
(2\pi)^{1/2}\left(\frac{T}{\mu_{\scriptscriptstyle \max}}\right)^{\!\!n/2}A_{n-1}\,-
\int^\infty_\rho\,e^{-\frac{\mu_{\scriptscriptstyle \max}}{\mathstrut 2T}\,r^{\scriptscriptstyle 2}}\,r^{n-1}\,\,dr
$$
where $A_n = E|X|^n$ for $X \sim  N(0,1)$, and that
$$
\int^\infty_\rho\,e^{-\frac{\mu_{\scriptscriptstyle \max}}{\mathstrut 2T}\,r^{\scriptscriptstyle 2}}\,r^{n-1}\,\,dr \leq 
\rho^{n-2}\,\frac{T}{\mu_{\scriptscriptstyle \max}}\,e^{-\frac{\mu_{\scriptscriptstyle \max}}{\mathstrut 2T}\,\rho^{\scriptscriptstyle 2}} \leq 
\rho^{n-2}\,\frac{T}{\mu_{\scriptscriptstyle \max}}\,e^{-\frac{U_{\scriptscriptstyle \rho}}{\mathstrut T}}
$$
Indeed, taken together, these give
$$
Z(T) \geq \omega_n\,\left(\frac{2}{\pi}\right)^{\!n-1}\,\left[(2\pi)^{1/2}\left(\frac{T}{\mu_{\scriptscriptstyle \max}}\right)^{\!\!n/2}A_{n-1}\,-
\rho^{n-2}\,\frac{T}{\mu_{\scriptscriptstyle \max}}\,e^{-\frac{U_{\scriptscriptstyle \rho}}{\mathstrut T}}\right]
$$
Finally, (\ref{eq:lower}) can be obtained by noting the second term in square brackets is negligeable compared to the first, as $T$ decreases to $0$, and by expressing $\omega_n$ and $A_{n-1}$ in terms of the Gamma function~\cite{wongbeals}. \hfill $\blacksquare$ \\[0.2cm] \end{subequations}
\indent \textbf{Proof of second estimate\,:} the Kantorovich distance between $P^\rho_{\scriptscriptstyle T}$ and the Dirac distribution $\delta_{\scriptscriptstyle x^*}$ is equal to the expectation of the distance to $x^*$, with respect to  $P^\rho_{\scriptscriptstyle T}$~\cite{villani}. Precisely,
\begin{subequations}
$$
W(P^\rho_{\scriptscriptstyle T},\delta_{\scriptscriptstyle x^*}) \,=\, \int_M d(x^*,z)\,P^\rho_{\scriptscriptstyle T}(dz)
$$
According to (\ref{eq:gibbs}) and (\ref{eq:prcon22}), this is
$$
W(P^\rho_{\scriptscriptstyle T},\delta_{\scriptscriptstyle x^*}) \,=\, \left(P_{\scriptscriptstyle T}(B_\rho)Z(T) \right)^{-1}\, \int_{B_\rho} d(x^*,z)\,e^{-\frac{U(z)}{\mathstrut T}}\,\mathrm{vol}(dz)
$$
Using (\ref{eq:gibbs}) to express the probability $P_{\scriptscriptstyle T}(B_\rho)\,$, this becomes
\begin{equation} \label{eq:prest21}
W(P^\rho_{\scriptscriptstyle T},\delta_{\scriptscriptstyle x^*}) \,=\, \frac{\int_{B_\rho} d(x^*,z)\,e^{-\frac{U(z)}{\mathstrut T}}\,\mathrm{vol}(dz)}{\int_{B_\rho}\,e^{-\frac{U(z)}{\mathstrut T}}\,\mathrm{vol}(dz)}
\end{equation}
A lower bound on the denominator can be found from (\ref{eq:later}) and subsequent inequalities, which were used to prove (\ref{eq:lower}). Precisely, these inequalities provide
\begin{equation} \label{eq:prest22}
 \int_{B_\rho}\,e^{-\frac{U(z)}{\mathstrut T}}\,\mathrm{vol}(dz) \,\geq\,\frac{1}{2}\,
\omega_n\,\left(\frac{2}{\pi}\right)^{\!n-1}\,(2\pi)^{1/2}\,A_{n-1}\left(\frac{T}{\mu_{\scriptscriptstyle \max}}\right)^{\!\!n/2}
\end{equation} 
whenever $T \leq T^1_o\,$. For the numerator in (\ref{eq:prest21}), it will be shown that, for any $T$,
\begin{equation} \label{eq:prest23}
  \int_{B_\rho} d(x^*,z)\,e^{-\frac{U(z)}{\mathstrut T}}\,\mathrm{vol}(dz)\,\leq\, \omega_n\,(2\pi)^{1/2}\,A_{n}\left(\frac{T}{\mu_{\scriptscriptstyle \min}}\right)^{\!\!(n+1)/2}
\end{equation}
Then, (\ref{eq:estimate2}) follows by dividing (\ref{eq:prest23}) by (\ref{eq:prest22}), and replacing in (\ref{eq:prest21}), after noting that
$A_n/A_{n-1} \,=\, \sqrt{2\pi}\,B^{-1}_n\,$. Thus, it only remains to prove (\ref{eq:prest23}). Using (\ref{eq:prcon21}), it is seen that
$$
  \int_{B_\rho} d(x^*,z)\,e^{-\frac{U(z)}{\mathstrut T}}\,\mathrm{vol}(dz)\,\leq\,   \int_{B_\rho} d(x^*,z)\,e^{-\frac{\mu_{\scriptscriptstyle \min}}{\mathstrut 2T}\,d^{\scriptscriptstyle 2}(x,x^*)}\,\mathrm{vol}(dz)
$$
By expressing this last integral in Riemannian spherical coordinates, as in (\ref{eq:comparison11}),
\begin{equation} \label{eq:prest24}
  \int_{B_\rho} d(x^*,z)\,e^{-\frac{U(z)}{\mathstrut T}}\,\mathrm{vol}(dz)\,\leq\, 
  \int^\rho_0\int_{S^{n-1}}\,r\,e^{-\frac{\mu_{\scriptscriptstyle \min}}{\mathstrut 2T}\,r^{\scriptscriptstyle 2}}\,\lambda(r,s)\,dr\,\omega_n(ds)
\end{equation}
From the volume comparison theorem in~\cite{chavel} (see Page 130), $\lambda(r,s) \leq r^{\!n-1}$. Therefore, (\ref{eq:prest24}) becomes
$$
  \int_{B_\rho} d(x^*,z)\,e^{-\frac{U(z)}{\mathstrut T}}\,\mathrm{vol}(dz)\,\leq\, \omega_n\,\int^\rho_0\,e^{-\frac{\mu_{\scriptscriptstyle \min}}{\mathstrut 2T}\,r^{\scriptscriptstyle 2}}\,r^n\,dr \leq \omega_n\,\int^\infty_0\,e^{-\frac{\mu_{\scriptscriptstyle \min}}{\mathstrut 2T}\,r^{\scriptscriptstyle 2}}\,r^n\,dr 
$$
The right-hand side is half the $n$th absolute moment of a normal distribution. Expressing this in terms of $A_n\,$, and replacing in (\ref{eq:prest24}), gives  (\ref{eq:prest23}). \hfill $\blacksquare$
\end{subequations}
\section{Proof of Lemma \ref{eq:lemma}}
Denote $G$ the connected component at identity of the group of isometries of $M$. It will be assumed that $G$ is simply-connected and semisimple~\cite{helgason}. Any geodesic $\gamma:I\rightarrow M$ is of the form~\cite{helgason}\cite{kobayashi},
\begin{subequations}
\begin{equation} \label{eq:geodesic}
  \gamma(t) \,=\, \exp(tY)\cdot x
\end{equation}
for some $x \in M$ and $Y \in \mathfrak{g}$, the Lie algebra of $G$, where $\exp:\mathfrak{g}\rightarrow G$ denotes the Lie group exponential mapping, and the dot denotes the action of $G$ on $M$. For each $t \in I$, the cut locus $\mathrm{Cut}(\gamma(t))$ of $\gamma(t)$ is given by
\begin{equation} \label{eq:cut1}
  \mathrm{Cut}(\gamma(t)) \,=\, \exp(tY)\cdot \mathrm{Cut}(x)
\end{equation}  
This is due to a more general result\,: let $M$ be a Riemannian manifold and $g:M\rightarrow M$ be an isometry of $M$. Then, $\mathrm{Cut}(g\cdot x) = g\cdot \mathrm{Cut}(x)$ for all $x \in M$. This is because $y \in \mathrm{Cut}(x)$ if and only if $y$ is conjugate to $x$ along some geodesic, or there exist two different geodesics connecting $x$ to $y$~\cite{chavel}\cite{kobayashi}. Both of these properties are preserved by the isometry $g$. 

In order to describe the set $\mathrm{Cut}(x)$, denote $K$ the isotropy group of $x$ in $G$, and $\mathfrak{k}$ the Lie algebra of $K$. Let $\mathfrak{g} = \mathfrak{k}+\mathfrak{p}$ be an orthogonal decomposition, with respect to the Killing form of $\mathfrak{g}$, and let $\mathfrak{a}$ be a maximal Abelian subspace of $\mathfrak{p}$. Define~$\mathcal{S} = K/C_{\mathfrak{a}}$ ($C_{\mathfrak{a}}$ the centraliser of $\mathfrak{a}$ in $K$), and consider the mapping
\begin{equation} \label{eq:phi}
  \phi(s,a) = \exp\left(\mathrm{Ad}(s)\,a\right)\cdot x \hspace{0.25cm}\mbox{ for }\, (s,a) \in \mathcal{S}\times \mathfrak{a}
\end{equation}
The set $\mathrm{Cut}(x)$ is the image under $\phi$ of a certain set $\mathcal{S}\times \partial Q\,$, which is now described, following~\cite{helgason}\cite{crittenden}.  

Let $\Delta_+$ be the set of positive restricted roots associated to the pair $(G,K)$, (each $\lambda \in \Delta_+$ is a linear form $\lambda:\mathfrak{a}\rightarrow \mathbb{R}$). Then, let $Q$ be the set of $a \in \mathfrak{a}$ such that $\left|\lambda(a)\right| \leq \pi$ for all $\lambda \in \Delta_+\,$, and $\partial Q$ the boundary of $Q$. Then 
\begin{equation} \label{eq:cutlocus}
  \mathrm{Cut}(x)\,=\, \phi(\mathcal{S}\times\partial Q)
\end{equation}
Recapitulating (\ref{eq:cut1}) and (\ref{eq:cutlocus}),
\begin{equation} \label{eq:cut2}
  \mathrm{Cut}(\gamma) \,=\, \Phi(I\times\mathcal{S}\times\partial Q) \hspace{0.1cm}\mbox{ where }\hspace{0.1cm} \Phi(t,s,a) = \exp(tY)\cdot \phi(s,a)
\end{equation}

Lemma \ref{eq:lemma} states that the topological dimension of $\mathrm{Cut}(\gamma)$ is strictly less than $\dim \, M$. This is proved using results from topological dimension theory~\cite{helgason}\cite{hurewic}.

Note that both $I$ and $\mathcal{S}$ are compact. Indeed, $\mathcal{S}$ is compact since it is the continuous image of the compact group $K$ under the projection $K \rightarrow K/C_{\mathfrak{a}}$. Also, $\partial Q$ is compact in $\mathfrak{a}$, and $\partial Q = \cup_{\lambda}\, \partial Q_\lambda$ where $\partial Q_\lambda= \partial Q  \cap \lbrace \lambda(a) = \pm\, \pi\rbrace$ for $\lambda \in \Delta_+\,$. Since $\lbrace \lambda(a) = \pm\,\pi\rbrace$ is the union of two (closed) hyperplanes~in~$\mathfrak{a}$, $\partial Q_\lambda$ is compact. Now, each $I \times \mathcal{S}\times \partial Q_\lambda$ is compact, and therefore closed. It follows from (\ref{eq:cut2}) that (see~\cite{hurewic}, Page 30),
\begin{equation} \label{eq:dim1}
  \dim\,\mathrm{Cut}(\gamma)\,=\, \dim\, \bigcup_\lambda\, \Phi(I\times\mathcal{S}\times\partial Q_\lambda)\leq\, \max_\lambda\, \dim\,  \Phi(I\times\mathcal{S}\times\partial Q_\lambda)
\end{equation}
But, for each $\lambda$, 
$$
\Phi(I\times\mathcal{S}\times\partial Q_\lambda) \,=\, \Phi(I\times\mathcal{S}_\lambda\times\partial Q_\lambda) \,\subset\,\Phi\left(\mathbb{R}\times\mathcal{S}_\lambda\times\lbrace \lambda(a) = \pm\,\pi\rbrace\right)
$$
where $\mathcal{S}_\lambda = K/C_{\lambda}$ ($C_{\lambda}$ the centraliser of $\lbrace \lambda(a) = \pm\,\pi\rbrace$ in $K$). The above inclusion implies (by~\cite{hurewic}, Page 26),
\begin{equation} \label{eq:dim2}
\dim\,\Phi(I\times\mathcal{S}\times\partial Q_\lambda) \,\leq \dim\,
\Phi\left(\mathbb{R}\times\mathcal{S}_\lambda\times \lbrace \lambda(a) = \pm\,\pi\rbrace\right)
\end{equation}
To conclude, note that the set $\mathbb{R} \times \mathcal{S}_\lambda \times \lbrace \lambda(a) = \pm\,\pi\rbrace$ is a differentiable manifold. It follows that (see~\cite{helgason}, Page 345),
\begin{equation} \label{eq:dim3}
  \dim \, \Phi\left(\mathbb{R} \times \mathcal{S}_\lambda \times\lbrace \lambda(a) = \pm\,\pi\rbrace\right) \, \leq \, \dim\left(\mathbb{R} \times \mathcal{S}_\lambda \times \lbrace \lambda(a) = \pm\,\pi\rbrace\right)
\end{equation}
The right-hand side of this inequality is
$$
\dim\left(\mathbb{R} \times \mathcal{S}_\lambda \times \lbrace \lambda(a) = \pm\,\pi\rbrace\right) = 1 + \dim\,\mathcal{S}_\lambda + \dim\,\mathfrak{a}-1
$$
since the dimension of a hyperplane in $\mathfrak{a}$ is $\dim\,\mathfrak{a}-1$. In addition, according to~\cite{helgason} (Page 296), $\dim\,\mathcal{S}_\lambda < \dim\,\mathcal{S}$. Thus,
$$
\dim\left(\mathbb{R} \times \mathcal{S}_\lambda \times \lbrace \lambda(a) = \pm\,\pi\rbrace\right) = \dim\,\mathcal{S}_\lambda + \dim\,\mathfrak{a} < 
\dim\,M
$$
since $\dim\,M = \dim\,\mathcal{S} + \dim\,\mathfrak{a}$~\cite{helgason}. Replacing this into (\ref{eq:dim3}), it follows from (\ref{eq:dim1}) and (\ref{eq:dim2}) that $\dim\,\mathrm{Cut}(\gamma) < \dim\,M$, as required.\hfill$\blacksquare$
\end{subequations}
\section{Proof of Corollary \ref{corr}}
The corollary can be split into the two following claims, which will be proved separately.\\[0.1cm]
\indent \textbf{First claim\,:} both integrals $G_x$ and $H_x$ converge for any value of $T$. \\[0.1cm]
\indent \textbf{Second claim\,:} $\mathcal{E}_{\scriptscriptstyle T}$ is $C^2$ throughout $M$, with derivatives given by (\ref{eq:derivatives}). \\[0.2cm]
The fact that $G_x$ and $H_x$ depend continuously on $x$ is contained in the second claim, since (\ref{eq:derivatives}) states that $G_x$ and $H_x$ are the gradient and Hessian of  $\mathcal{E}_{\scriptscriptstyle T}$ at $x$. \\[0.2cm]
\indent In the following proofs, the notation $\mathrm{D}(x) = M - \mathrm{Cut}(x)$ will be used, in~order to avoid cumbersome expressions.

\vfill\pagebreak
\indent \textbf{Proof of first claim\,:}  The convergence of the integral $G_x$ is straightforward, since the integrand $G_x(z)$ is a smooth and bounded function, from $\mathrm{D}(x)$ to $T_xM$. This is because, by definition, $G_x(z)$ is given by
\begin{equation} \label{eq:expoponential}
G_x(z) \,=\, -\,\mathrm{Exp}^{-1}_x(z)
\end{equation}
\begin{subequations}
where $\mathrm{Exp}$ is the Riemannian exponential mapping~\cite{chavel}. Therefore, $G_x(z)$ is smooth. In addition, $G_x(z)$ is bounded, in Riemannian norm, by $\mathrm{diam}\, M$. 

The convergence of the integral $H_x$ is more difficult. While the integrand $H_x(z)$ is smooth on $\mathrm{D}(x)$, it is not bounded. It will be seen that $H_x$ is an absolutely convergent improper integral. 

Recall the mapping $\phi$ defined in (\ref{eq:phi}). Let $D_+$ be the set of points $a \in \mathfrak{a}$ which belong to the interior of $Q$, and which verify $\lambda(a) \geq 0$ for each $\lambda \in \Delta_+\,$.\hfill\linebreak Let ${D}^{\scriptscriptstyle o}_+$ be the interior of $D_+\,$. Then, $\phi$ maps $\mathcal{S}\times D_+$ onto $\mathrm{D}(x)$, and is a diffeomorphism of $\mathcal{S}\times {D}^{\scriptscriptstyle o}_+$ onto its image in $\mathrm{D}(x)$~\cite{helgason}\cite{crittenden} (see Chapter VII in~\cite{helgason}). Using Sard's theorem~\cite{bogachev}, it follows from the definition of $H_x$ that
\begin{equation} \label{eq:prcor1}
  H_x \,=\, \int_{\mathcal{S}}\int_{{D}_+} H_x(\phi(s,a))\,p_{\scriptscriptstyle T}(\phi(s,a))\,J(a)\,da\,\omega(ds)
\end{equation}
where $p_{\scriptscriptstyle T}$ denotes the density of $P_{\scriptscriptstyle T}$ with respect to the Riemannian volume of~$M$, and $J(a)$ is the Jacobian determinant of $\phi\,$, given by~\cite{helgason}
\begin{equation} \label{eq:ja}
J(a) \,=\, \prod_{\lambda \in \Delta_+}\,\left(\sin\lambda(a)\right)^{m_\lambda}
\end{equation}
with $m_\lambda$ the multiplicity of the restricted root $\lambda$, and where $\omega(ds)$ is the invariant Riemannian volume induced on $\mathcal{S}$ from $K$. 

Now, $H_x(\phi(s,a))$ can be expressed as follows ($\cot$ is the cotangent function)
\begin{equation} \label{eq:hessa}
 H_x(\phi(s,a)) \,=\, \Pi_0(s) \,+\, \sum_{\lambda \in \Delta_+}\,\lambda(a)\cot\lambda(a)\, \Pi_{\lambda}(s)
\end{equation}
where $\Pi_0(s)$ and the $\Pi_{\lambda}(s)$ denote orthogonal projectors, onto the respective eigenspaces of $H_x(\phi(s,a))$. 

According to this expression, $H_x(\phi(s,a))$ diverges to $-\infty$ whenever $\lambda(a) = \pi$. However, the product 
$$
H_x(\phi(s,a))\,p_{\scriptscriptstyle T}(\phi(s,a))\,J(a)
$$ 
which appears under the integral in (\ref{eq:prcor1}), is clearly continuous and bounded on the domain of integration. Thus, the absolute convergence of the integral $H_x$ follows immediately from (\ref{eq:prcor1}). 
\end{subequations}
It now remains to provide a proof of (\ref{eq:hessa}). This is here only briefly indicated. Expression (\ref{eq:hessa}) is a slight improvement of the one in~\cite{joao} (see Theorem IV.1, Page 636), where it is enough to note that if $R$ is the curvature tensor of $M$, then the operator $R_v(u) = R(v,u)v$ has the  eigenvalues $0$ and $(\lambda(a))^2$ for each $\lambda \in \Delta_{+\,}$, whenever $v,u \in T_xM \simeq \mathfrak{p}$ with $v = \mathrm{Ad}(s)\,a$~\cite{helgason}\cite{crittenden}. It is well-known, by properties of the Jacobi equation~\cite{chavel}, that $H_x(\phi(s,a))$ has the same eigenspace decomposition as $R_v$, in this case. \hfill$\blacksquare$ \vfill\pagebreak
\indent \textbf{Proof of second claim\,:} the proof of this claim relies in a crucial way on Lemma \ref{eq:lemma}. To compute the gradient and Hessian of the function $\mathcal{E}_{\scriptscriptstyle T}$ at $x \in M$, consider any geodesic $\gamma:I\rightarrow M$, defined on a compact interval $I = [-\tau,\tau]$, such that $\gamma(0) = x$. For each $t \in I$, by definition of the function $\mathcal{E}_{\scriptscriptstyle T\,}$,
\begin{subequations}
\begin{equation} \label{eq:prcor21}
  \mathcal{E}_{\scriptscriptstyle T}(\gamma(t)) \,=\, \frac{1}{2}\,\int_M\,d^2(\gamma(t),z)\,P_{\scriptscriptstyle T}(dz)
\end{equation}
However, Lemma \ref{eq:lemma} states that the set
$$
\mathrm{Cut}(\gamma) \,=\, \bigcup_{t\in I}\, \mathrm{Cut}(\gamma(t))
$$
has Riemannian volume equal to zero. From (\ref{eq:gibbs}), it is clear that $P_{\scriptscriptstyle T}$ is absolutely continuous with respect to Riemannian volume. Therefor, $\mathrm{Cut}(\gamma)$ can be removed from the domain of integration in (\ref{eq:prcor21}). Then,
\begin{equation} \label{eq:prcor22}
 \mathcal{E}_{\scriptscriptstyle T}(\gamma(t)) \,=\, \frac{1}{2}\,\int_{\mathrm{D}(\gamma)}\,d^2(\gamma(t),z)\,P_{\scriptscriptstyle T}(dz)
\end{equation}
where $\mathrm{D}(\gamma) = M - \mathrm{Cut}(\gamma)$. Now, for each $z \in \mathrm{D}(\gamma)$, the function 
$$
t\mapsto f_z(t) = \frac{1}{2}\,d^2(\gamma(t),z)
$$ 
is twice continuously differentiable with respect to $t \in I$, with
\begin{equation} \label{eq:prcor23}
  \frac{df_z}{dt}\,=\left\langle G_{\gamma(t)}(z),\dot{\gamma}\right\rangle \hspace{0.15cm}\mbox{ and }\hspace{0.15cm}
 \frac{d^2f_z}{dt^2}\,=H_{\gamma(t)}(z)\left(\dot{\gamma},\dot{\gamma}\right)
\end{equation}
where $\langle\cdot,\cdot\rangle$ denotes the Riemannian metric of $M$, and $\dot{\gamma}$ the velocity of the geodesic $\gamma$. Indeed, this holds because the geodesic $\gamma$ does not intersect the cut locus $\mathrm{Cut}(z)$~(see \cite{chavel}).  \\[0.1cm]
\indent The claim that $\mathcal{E}_{\scriptscriptstyle T}$ is twice differentiable, and has derivatives given by (\ref{eq:derivatives}), follows from (\ref{eq:prcor22}) and (\ref{eq:prcor23}), by differentiation under the integral sign, provided it can be shown that the families of functions
$$
\left \lbrace\,z\mapsto G_{\gamma(t)}(z)\,;\,t\in I\,\right\rbrace \hspace{0.15cm}\mbox{ and }\hspace{0.15cm}
\left \lbrace\,z\mapsto H_{\gamma(t)}(z)\,;\,t\in I\,\right\rbrace
$$
which all have the common domain of definition $\mathrm{D}(\gamma)$, are uniformly integrable with respect to $P_{\scriptscriptstyle T}$~\cite{bogachev}. Roughly, uniform integrability means that the rate of absolute convergence of the following integrals does not depend on $t$, 
$$
G_{\gamma(t)}\,=\,\int_{\mathrm{D}(\gamma)}\,G_{\gamma(t)}(z)\,P_{\scriptscriptstyle T}(dz) \hspace{0.25cm};\hspace{0.25cm}
H_{\gamma(t)}\,=\,\int_{\mathrm{D}(\gamma)}\,H_{\gamma(t)}(z)\,P_{\scriptscriptstyle T}(dz)
$$ 
This is clear for the integrals $G_{\gamma(t)}$ because $G_{\gamma(t)}(z)$ is bounded in Riemannian norm by $\mathrm{diam}\,M$, uniformly in $t$ and $z$ (see the proof of the first claim). 

Then, consider the integral $H_x = H_{\gamma(0)}$, and recall Formulae (\ref{eq:prcor1}) and (\ref{eq:hessa}). Each $z \in \mathrm{D}(\gamma)$ can be written under the form $z = \phi(s,a)$ where $(s,a) \in \mathcal{S}\times D_{+\,}$. Accordingly, it follows from (\ref{eq:hessa}) that
\begin{equation} \label{eq:frobenius}
 \left \Vert H_x(z) \right\Vert_{F} \,\leq (\dim\,M)^{\frac{1}{2}}\, \max\left\lbrace 1, \left|\kappa(a)\cot\kappa(a)\right|\right\rbrace
\end{equation}
where $\Vert\cdot\Vert_F$ is the Frobenius norm with respect to the Riemannian metric of $M$, and $\kappa \in \Delta_+$ is the highest restricted root~\cite{helgason} ($\kappa(a) \geq \lambda(a)$ for $\lambda \in \Delta_{+\,}$, $a \in D_+$). 

The required uniform integrability is equivalent to the statement that 
\begin{equation} \label{eq:prcor24}
\lim_{K\rightarrow \infty}\,\,  \int_{\mathrm{D}(\gamma)}\,  \left \Vert H_x(z) \right\Vert_{F}\,\mathbf{1}\left\lbrace \left \Vert H_x(z) \right\Vert_{F} > K \right\rbrace\,P_{\scriptscriptstyle T}(dz) \,=\, 0
\end{equation}
where the rate of convergence to this limit does not depend on $x$. But, according to (\ref{eq:frobenius}), if $K > 1$, there exists $\epsilon > 0$ such that
$$
\left\lbrace \left \Vert H_x(z) \right\Vert_{F} > K \right\rbrace \,=\,
\left\lbrace \kappa(a) > \pi-\epsilon \right\rbrace
$$ 
and $\epsilon \rightarrow 0$ as $K \rightarrow \infty$. In this case, the integral in (\ref{eq:prcor24}) is less than
\begin{equation} \label{eq:prcor25}
(\dim\,M)^{\frac{1}{2}}\,\left(\sup_z\,p_{\scriptscriptstyle T}(z)\right)\,\int_{\mathcal{D}(\gamma)}\,\left|\kappa(a)\cot\kappa(a)\right|\,\mathbf{1}\left\lbrace \kappa(a) > \pi-\epsilon \right\rbrace\,\mathrm{vol}(dz)
\end{equation}
Now, using the same integral formula as in (\ref{eq:prcor1}), this last integral is equal to
$$
\begin{array}{l}
\int_{\mathcal{S}}\,\int_{D_+}\,\left|\kappa(a)\cot\kappa(a)\right|\,\mathbf{1}\left\lbrace \kappa(a) > \pi-\epsilon \right\rbrace\,J(a)\,da\,\omega(ds) \,=\, 
\\[0.3cm]
\omega(\mathcal{S})\,\int_{D_+}\,\left[\,\left|\kappa(a)\cot\kappa(a)\right|\,J(a)\,\right]\,\mathbf{1}\left\lbrace \kappa(a) > \pi-\epsilon \right\rbrace\,da
\end{array}
$$
In view of (\ref{eq:ja}), since $\kappa \in \Delta_{+\,}$, the function in square brackets is bounded on the closure of $D_+$. In fact~\cite{helgason}, its supremum is $\kappa^2
= (\kappa,\kappa)$ where $(\cdot,\cdot)$ is the scalar product induced on $\mathfrak{a}^*$ (the dual space of $\mathfrak{a})$ by the Killing form of $\mathfrak{g}$. Finally, by (\ref{eq:prcor25}), the integral in (\ref{eq:prcor24}) is less than
$$
(\dim\,M)^{\frac{1}{2}}\,\left(\sup_z\,p_{\scriptscriptstyle T}(z)\right)\,\omega(\mathcal{S})\,\kappa^2\,
\int_{D_+}\,\mathbf{1}\left\lbrace \kappa(a) > \pi-\epsilon \right\rbrace\,da 
$$
Since, $\kappa(a) \in [0,\pi)$ for $a \in D_+\,$, this last integral converges to $0$ as  $\epsilon \rightarrow 0$, at a rate which does not depend on $x$. This proves the required uniform integrability, so the proof is now complete. \hfill$\blacksquare$ 
\end{subequations}
\section{Proof of Proposition \ref{prop:convexity}}
\begin{subequations}

\indent \textbf{Remark\,:} in the statement of Proposition \ref{prop:convexity}, the notation $\kappa^2$ is used for the maximum sectional curvature of $M$. In the previous proof of Corollary \ref{corr}, the same notation $\kappa^2$ was used for the squared norm of the highest restricted root. This is not an abuse of notation, since the two quantities are in fact equal~\cite{helgason} (see Page 334). \\[0.2cm]
\textbf{Proof of \textit{(i)}\,:} let $x \in B(x^*,\delta)$. By (\ref{eq:derivatives}) of Corollary \ref{corr}, $\nabla^2 \mathcal{E}_{\scriptscriptstyle T}(x)$ is equal to $H_x\,$. To obtain (\ref{eq:convexity1}), decompose $H_x$ into two integrals
\begin{equation} \label{eq:prconvex1}
  H_x  \,=\, \int_{B(x,r_{cx})}\, H_x(z)\,P_{\scriptscriptstyle T}(dz) \,+\,\int_{\mathrm{D}(x)-B(x,r_{cx})}\, H_x(z)\,P_{\scriptscriptstyle T}(dz)
\end{equation}
This is possible since $B(x,r_{cx}) \subset \mathrm{D}(x)$, where $\mathrm{D}(x) = M - \mathrm{Cut}(x)$. The first integral in (\ref{eq:prconvex1}) will be denoted $I_{\scriptscriptstyle 1\,}$, and the second integral $I_{\scriptscriptstyle 2\,}$. 

With regard to $I_{\scriptscriptstyle 1\,}$, note the inclusions $B(x^*,\delta) \subset B(x,2\delta) \subset B(x,r_{cx})$, which follow from the triangle inequality. In addition, note that $H_x(z) \geq 0$ (in~the~Loewner order~\cite{matrix}), for $z \in B(x,r_{cx})$. Therefore,
\begin{equation} \label{eq:prconvex2}
  I_{\scriptscriptstyle 1}\,\geq\, \int_{B(x^*,\delta)}\,H_x(z)\,P_{\scriptscriptstyle T}(dz)
\end{equation}
However, from (\ref{eq:hessa}) and the definition of $\kappa \in \Delta_+\,$, 
\begin{equation} \label{eq:prconvex3}
H_x(z) \geq \kappa(a)\cot\kappa(a)
\end{equation}
for $z = \phi(s,a) \in \mathrm{D}(x)$. Using the Cauchy-Scwharz inequality, $\kappa(a) \leq \kappa\, \Vert a \Vert$. Moreover, (\ref{eq:phi}) implies $\Vert a\Vert = d(x,z)$, since $\mathrm{Ad}(s)$ is an isometry. Accordingly, if $z \in B(x,2\delta)$, it follows from (\ref{eq:prconvex3})
\begin{equation} \label{eq:prconvex4}
H_x(z) \geq \kappa(a)\cot\kappa(a) \geq 2\kappa\delta\cot(2\kappa\delta) \,=\, \mathrm{Ct}(2\delta)\, > 0
\end{equation} 
where the last inequality is because $2\delta < r_{cx} = \kappa^{-1}\frac{\pi}{2\,}$. Replacing in (\ref{eq:prconvex2}) gives
$$
I_{\scriptscriptstyle 1}\,\geq\, \mathrm{Ct}(2\delta)\,P_{\scriptscriptstyle T}(B(x^*,\delta)) \,=\, \mathrm{Ct}(2\delta)\,\left[1 - P_{\scriptscriptstyle T}(B^c(x^*,\delta))\right]
$$
Finally, (\ref{eq:prest12}) and (\ref{eq:lower}) imply that $P_{\scriptscriptstyle T}(B^c(x^*,\delta)) \leq \mathrm{vol}(M)\,f(T)$, where $f(T)$ was defined in (\ref{eq:fnt}) -- Precisely, this follows after replacing $\rho$ by $\delta$ in (\ref{eq:prest12}). Thus,
\begin{equation} \label{eq:prconvex5}
I_{\scriptscriptstyle 1}\,\geq\, \mathrm{Ct}(2\delta)\left( 1 - \mathrm{vol}(M) f(T)\right)
\end{equation}
\end{subequations}
The proof of (\ref{eq:convexity1}) will be completed by showing
\begin{subequations}
\begin{equation} \label{eq:pri21}
 I_{\scriptscriptstyle 2}\,\geq\, -\, \pi  A_M f(T)
\end{equation}
To show this, note using (\ref{eq:prconvex3}) that
\begin{equation} \label{eq:pri22}
I_{\scriptscriptstyle 2}\,\geq\, \int_{\mathrm{D}(x)-B(x,r_{cx})}\,\kappa(a)\cot\kappa(a)\,P_{\scriptscriptstyle T}(dz)
\end{equation}
Now, $\kappa(\alpha)\cot\kappa(\alpha)$ is negative if and only if $\kappa(\alpha) \geq \frac{\pi}{2\,}$. However, the set of $z  = \phi(s,a)$ where $\kappa(a) \geq \frac{\pi}{2\,}$ is a subset of $\mathrm{D}(x)-B(x,r_{cx})$. Indeed, $\kappa(a) \geq \frac{\pi}{2}$ implies $\Vert a \Vert \geq \kappa^{-1}\frac{\pi}{2} = r_{cx\,}$, by the Cauchy-Schwarz inequality, and this is the same as $d(x,z) \geq r_{cx\,}$, since $\Vert a\Vert = d(x,z)$. Therefore, it follows from (\ref{eq:pri22}),
\begin{equation} \label{eq:pri33}
I_{\scriptscriptstyle 2}\,\geq\, \int_{\mathrm{D}(x)}\,\mathbf{1}\lbrace \kappa(a) \geq \pi/2\rbrace\,\kappa(a)\cot\kappa(a)\,P_{\scriptscriptstyle T}(dz)
\end{equation}
Using the same integral formula as in (\ref{eq:prcor1}), this last integral is equal to
$$
\begin{array}{l}
\int_{\mathcal{S}}\,\int_{D_+}\,\,\mathbf{1}\lbrace \kappa(a) \geq \pi/2\rbrace\,\,\kappa(a)\cot\kappa(a)\,p_{\scriptscriptstyle T}(\phi(s,a))\,J(a)\,da\,\omega(ds) \geq \\[0.2cm]
-\, \int_{\mathcal{S}}\,\int_{D_+}\,\,\mathbf{1}\lbrace \kappa(a) \geq \pi/2\rbrace\,\kappa(a)\,p_{\scriptscriptstyle T}(\phi(s,a))\,da\,\omega(ds)
\end{array}
$$
because the product $\cot\kappa(a)\,J(a) \geq -1$ for all $a \in D_{+\,}$. Using this last inequality, and the fact that $\kappa(a) \leq \pi$ for all $a \in D_{+\,}$, it follows from (\ref{eq:pri33}),
\begin{equation} \label{eq:pri34}
I_{\scriptscriptstyle 2}\,\geq\, -\pi\,\int_{\mathcal{S}}\,\int_{D_+}\,\mathbf{1}\lbrace \kappa(a) \geq \pi/2\rbrace\,p_{\scriptscriptstyle T}(\phi(s,a))\,da\,\omega(ds)
\end{equation}
Recall that $\lbrace \kappa(a) \geq \pi/2\rbrace \subset B^c(x,r_{cx})$, as discussed before (\ref{eq:pri33}). In particular, this implies $\lbrace \kappa(a) \geq \pi/2\rbrace \subset B^c(x^*,\delta)$. However, by (\ref{eq:gibbs}) and (\ref{eq:lower}), $p_{\scriptscriptstyle T}(z) \leq f(T)$ for all $z \in B^c(x^*,\delta)$. Returning to (\ref{eq:pri34}), this gives
\begin{equation} \label{eq:pri35}
I_{\scriptscriptstyle 2}\,\geq\, -\pi\,f(T)\,\int_{\mathcal{S}}\,\int_{D_+}\,da\,\omega(ds)
\end{equation}
The double integral on the right-hand side is a constant which depends only on the structure of the symmetric space $M$. Denoting this constant by $A_M$ gives the required lower bound (\ref{eq:pri21}), and completes the proof of (\ref{eq:convexity1}). \hfill$\blacksquare$ \end{subequations}\\[0.2cm]
\textbf{Proof of \textit{(ii)}\,:} fix $\delta < \frac{1}{2}r_{cx\,}$, and let $T_{\scriptscriptstyle \delta}$ be given by (\ref{eq:camille}). If $T \leq T_{\scriptscriptstyle \delta\,}$, then $T < T^2_{\scriptscriptstyle \delta\,}$, so the definition of $T^2_{\scriptscriptstyle \delta}$ implies
\begin{subequations}
\begin{equation} \label{eq:pr221}
f(T) \,<\, \frac{\mathrm{Ct}(2\delta)}{\mathrm{Ct}(2\delta)\,\mathrm{vol}(M) + \pi A_M}
\end{equation}
Now, by (\ref{eq:convexity1}), 
\begin{equation} \label{eq:pr222}
\nabla^2 \mathcal{E}_{\scriptscriptstyle T}(x) \geq  \mathrm{Ct}(2\delta)\left( 1 - \mathrm{vol}(M) f(T)\right) - \pi  A_M f(T) 
\end{equation}
for all $x \in B(x^*,\delta)$. However, it is clear from (\ref{eq:pr221}), that the right-hand side of this inequality is strictly positive. It follows that $\mathcal{E}_{\scriptscriptstyle{T}}$ is strongly convex on $B(x^*,\delta)\,$. Thus, to complete the proof, it only remains to show that any global minimum $\bar{x}_{\scriptscriptstyle T}$ of $\mathcal{E}_{\scriptscriptstyle{T}}$ must belong to $B(x^*,\delta)$. Indeed, since $\mathcal{E}_{\scriptscriptstyle{T}}$ is strongly convex on $B(x^*,\delta)\,$, it has only one local minimum in $B(x^*,\delta)$. Therefore, $\mathcal{E}_{\scriptscriptstyle{T}}$ can have only one global minimum $\bar{x}_{\scriptscriptstyle T\,}$. 

By \textit{(i)} of Proposition \ref{prop:concentration}, to prove that $\bar{x}_{\scriptscriptstyle T} \in B(x^*,\delta)$, it is enough to prove 
\begin{equation} \label{eq:pr223}
W(P_{\scriptscriptstyle{T\,}},\delta_{\scriptscriptstyle x^*}) < \frac{\delta^2}{(4\,\mathrm{diam}\, M)}
\end{equation}
However, if $T \leq T_{\scriptscriptstyle \delta\,}$, then $T < T_o\,$. Therefore, by \textit{(ii)} of Proposition \ref{prop:concentration}, $W(P_{\scriptscriptstyle{T\,}},\delta_{\scriptscriptstyle x^*})$ satisfies inequality (\ref{eq:concentration2}). Furthermore, because $T < T^1_{\scriptscriptstyle \delta\,}$, it follows from the definition of $T^1_{\scriptscriptstyle \delta}$ that
$$
\sqrt{2\pi}\,(T/\mu_{\scriptscriptstyle\min})^{1/2}\,<\,\delta^2\,\left(\mu_{\scriptscriptstyle\min}/\mu_{\scriptscriptstyle\max}\right)^{n/2}\,D_n
$$
or, by replacing the expression of $D_n\,$, and simplifying
\begin{equation} \label{eq:pr224}
\sqrt{2\pi}\,\left(\pi/2\right)^{n-1}\,B^{-1}_n\,\left(\mu_{\scriptscriptstyle \max}/\mu_{\scriptscriptstyle \min}\right)^{n/2}\,\left(T/\mu_{\scriptscriptstyle \min}\right)^{1/2}\,\,  <\,\, \frac{\delta^2}{(4\,\mathrm{diam}\, M)}
\end{equation}
Thus, (\ref{eq:pr223}) follows from (\ref{eq:concentration2}) and (\ref{eq:pr224}). This proves that $\bar{x}_{\scriptscriptstyle T}$ belongs to $B(x^*,\delta)$, and therefore $\bar{x}_{\scriptscriptstyle T}$ is the unique global minimum of $\mathcal{E}_{\scriptscriptstyle T\,}$. But this is equivalent to saying that 
$\bar{x}_{\scriptscriptstyle T}$ is the unique barycentre of $P_{\scriptscriptstyle T\,}$. \hfill$\blacksquare$ \end{subequations}
\section{Proof of Proposition \ref{prop:main}}
fix $\delta < \frac{1}{2}r_{cx\,}$, and let $T_{\scriptscriptstyle \delta}$ be given by (\ref{eq:camille}). By \textit{(ii)} of Proposition \ref{prop:convexity}, if $T \leq T_{\scriptscriptstyle \delta\,}$, then $\mathcal{E}_{\scriptscriptstyle T}$ is strictly convex on $B(x^*,\delta)$, with unique global minimum $\bar{x}_{\scriptscriptstyle T} \in B(x^*,\delta)\,$. By definition, this unique global minimum $\bar{x}_{\scriptscriptstyle T}$ is the unique barycentre of $P_{\scriptscriptstyle T\,}$. 

Accordingly, to prove that $\bar{x}_{\scriptscriptstyle T} = x^*$, it is enough to prove that $x^*$ is a stationary point of $\mathcal{E}_{\scriptscriptstyle T\,}$. Indeed,  as $\mathcal{E}_{\scriptscriptstyle T}$ is strictly convex on $B(x^*,\delta)$, it can have only one stationary point in $B(x^*,\delta)\,$. This stationary point is then identical to $\bar{x}_{\scriptscriptstyle T\,}$. 

The fact that $x^*$ is a stationary point of $\mathcal{E}_{\scriptscriptstyle T}$ will follow because $U$ is invariant by geodesic symmetry about $x^*$. This invariance will be seen to imply
\begin{subequations}
\begin{equation} \label{eq:prprop31}
  ds_{\scriptscriptstyle x^*}\cdot\, G_{x^*} \,=\, G_{x^*}
\end{equation}
which is equivalent to $G_{x^*} = 0$, since the derivative $ds_{\scriptscriptstyle x^*}$ is equal to minus the identity, on the tangent space $T_{\scriptscriptstyle x^*}M$~\cite{helgason}. By (\ref{eq:derivatives}) of Corollary \ref{corr}, this shows that $\nabla \mathcal{E}_{\scriptscriptstyle T}(x^*) = 0$, so $x^*$ is indeed a stationary point of $\mathcal{E}_{\scriptscriptstyle T\,}$.

To obtain (\ref{eq:prprop31}), it is possible to write, from the definition of $G_{x^*}$,
\begin{equation} \label{eq:prprop32}
ds_{\scriptscriptstyle x^*}\cdot\, G_{x^*} \,=\, ds_{\scriptscriptstyle x^*}\cdot\, 
\int_{\mathrm{D}(x)}\,G_{x^*}(z)\,P_{\scriptscriptstyle T}(dz)
\end{equation}
where $\mathrm{D}(x) = M - \mathrm{Cut}(x)$. From (\ref{eq:expoponential}), since $s_{\scriptscriptstyle x^*}$ is an isometry, and reverses geodesics passing through $x^*$,
$$
ds_{\scriptscriptstyle x^*}\cdot\, G_{x^*}(z)\,=\, G_{x^*}(s_{\scriptscriptstyle x^*}(z)) 
$$
Replacing this into (\ref{eq:prprop32}), and using $w = s_{\scriptscriptstyle x^*}(z)$ as a new variable of integration, it follows that
\begin{equation} \label{eq:prprop33}
ds_{\scriptscriptstyle x^*}\cdot\, G_{x^*}
\,=\, \int_{\mathrm{D}(x)}\,G_{x^*}(w)\,\left(P_{\scriptscriptstyle T}\circ s_{\scriptscriptstyle x^*}\right)(dw)
\end{equation}
because $s^{\scriptscriptstyle -1}_{\scriptscriptstyle x^*} = s_{\scriptscriptstyle x^*}$ and $s_{\scriptscriptstyle x^*}$ maps $\mathrm{D}(x)$ onto itself. 
Now, note that $P_{\scriptscriptstyle T}\circ s_{\scriptscriptstyle x^*} = P_{\scriptscriptstyle T\,}$. This is clear, since from (\ref{eq:gibbs}),
$$
\left(P_{\scriptscriptstyle T}\circ s_{\scriptscriptstyle x^*}\right)(dw) \,=\, 
\left(Z(T)\right)^{-1}\,\exp\left[-\frac{\left(U \circ s_{\scriptscriptstyle x^*}\right)\! (w)}{T}\right]\left(\mathrm{vol}\circ s_{\scriptscriptstyle x^*}\right)(dw)
$$
However, by assumption, $U \circ s_{\scriptscriptstyle x^*} (w) = U(w)$. Moreover, since $s_{\scriptscriptstyle x^*}$ is an isometry, it preserves Riemannian volume, so $\left(\mathrm{vol}\circ s_{\scriptscriptstyle x^*}\right)(dw) = \mathrm{vol}(dw)$. Thus, (\ref{eq:prprop33}) reads
$$
ds_{\scriptscriptstyle x^*}\cdot\, G_{x^*}
\,=\, \int_{\mathrm{D}(x)}\,G_{x^*}(w)\,P_{\scriptscriptstyle T}(dw)
$$
By definition, the right-hand side is $G_{x^*}$, so (\ref{eq:prprop31}) is obtained.\hfill$\blacksquare$

\end{subequations}

\bibliographystyle{splncs04}

\end{document}